\documentclass{elsart}
\theoremstyle{definition}
\newtheorem{theorem}{Theorem}[section]

\newtheorem{corollary}{Corollary}[section]

\theoremstyle{remark}
\usepackage[mathscr]{eucal}
\usepackage{mathrsfs}
\usepackage{layout}
\usepackage{epstopdf}
\usepackage{graphicx}
\usepackage{float}
\usepackage{amssymb}
\begin{document}
\begin{frontmatter}
\title{RESULTS ON THE SUPREMUM OF FRACTIONAL BROWNIAN MOTION}
\author{CEREN VARDAR}\footnote{This study is supported by Tubitak with project number
110T674}
\address{Department of Mathematics,\
TOBB Economy and Technology University, Ankara, Turkey}
\begin{abstract}
We show that the distribution of the square of the supremum of
reflected fractional Brownian motion up to time a, with Hurst
parameter-H greater than 1/2, is related to the distribution of its
hitting time to level $1,$ using the self similarity property of
fractional Brownian motion. It is also proven that second moment of
supremum of reflected fractional Brownian motion up to time $a$ is
bounded above by $a^{2H}.$ Similar relations are obtained for the
supremum of fractional Brownian motion with Hurst parameter greater
than 1/2, and its hitting time to level $1.$ What is more, we obtain
an upper bound on the complementary probability distribution of the
supremum of fractional Brownian motion and reflected fractional
Brownian motion up to time a, using Jensen's and Markov's
inequalities. A sharper bound is observed on the distribution of the
supremum of fractional Brownian motion by the properties of Gamma
distribution. Finally, applications of the given results to the
financial markets are investigated and partial results are provided.
\end{abstract}
\end{frontmatter}


\section{Introduction}
\subsection*{Stochastic Model}
The price of the stock, a volatile asset, is often described by
Black-Scholes model which is also known as Geometric Brownian model.
\newpage
The stochastic differential equation of this model is
\begin{eqnarray}dY_{t}=\bar{\mu}Y_{t} dt+ \bar{\sigma}
Y_{t}dW_{t},~~~~~~0\leq t \leq T\end{eqnarray} where $\bar{\mu}$ is
a constant mean rate of return and $\bar{\sigma}$ is a constant
volatility, $W$ is the standard Brownian motion and $T$ is the time
of maturity of the stock. When we solve this equation explicitly, we
obtain
\begin{eqnarray} Y_{t}= Y_{0}\exp{((\bar{\mu}-\frac{1}{2}\bar{\sigma}^{2})t+\bar{\sigma} W_{t})}\end{eqnarray}
where $Y_{0}$ is the initial value of one share of stock. However,
this model has well known deficiencies such as increments of
standard Brownian motion are independent. This is not what we are
exposed to in real life. In order to overcome this problem various
alternative models have been suggested such as fractional Brownian
motion, (fBm).

Fractional Brownian motion also appears naturally in many
situations. Some examples are, the level of water in a river as a
function of time, the characters of solar activity as a function of
time, the values of the log returns of a stock, the prices of
electricity in a liberated electricity market, \cite{bib:oksendal}.
In the first three examples fBm with Hurst parameter, $H>1/2$ is
used which means that the process is persistent. And the last
example is modeled by fBm with Hurst parameter $H<1/2.$ Through out
this article we are interested in the financial applications of fBm,
specifically the model of the values of log returns of a stock. In
other words we are interested in the fBm with Hurst parameter
$H>1/2.$ The Black-Scholes model for the values of the log returns
of a stock using fBm is given as
\begin{eqnarray} Y_{t}= Y_{0}\exp{((r+\mu)t+\sigma B^{H}_{t})},~~~~~~  0\leq t \leq T\end{eqnarray}
where $Y_{0}$ is the initial value, $r$ is the constant interest
rate, $\mu$ is the constant drift and $\sigma$ is the constant
diffusion coefficient of fBm which is denoted by
$(B^{H}_{t})_{t\geq0}.$

The advantage of modeling with fBm to the other models is its
capability of displaying the dependence between returns on different
days. FBm on the other hand is not a semimartingale and is not a
Markov process and it allows arbitrage, \cite{bib:rogers}. However,
in order to display long range dependence and to have semimartingale
property other models can be constructed using fBm as a central
model, \cite{bib:rogers}. Due to the complex nature of fBm, the most
useful and efficient classical mathematical techniques for
stochastic calculus are not available for fBm. Therefore most of the
results in the literature are given as bounds on the characteristics
of fBm. In spite of this, fBm still possess some nice properties, it
is a Gaussian process and it holds self similarity property. Most of
the techniques developed for fBm are related to these properties.
Now, there is no doubt that the investors would be interested in
having some information on value of the supremum of the asset in
order to manage the risk or maybe to hedge financial assets and to
construct portfolios. Therefore in this study, mainly two new
results on the supremum of fBm and an application of these results
are introduced.
\begin{itemize}
  \item An identity for the distribution of square of the supremum
  of fBm
  \item An identity for the distribution of square of the supremum
  of  reflected fBm
  \item An upper bound on the second moment of reflected fBm
  \item A lower bound on the distribution of the supremum of fBm upto time
  $1$ or any fixed time $t.$
  \item A theoretical application of the above results is given, which
is commonly used in the literature to obtain a measure of risk in
finance that is called as "maximum drawdown."
\end{itemize}

\subsection*{Fractional Brownian Motion}
Let us start with introducing fBm. FBm was first introduced within
Hilbert space framework by Kolmogorov, \cite{bib:kolmogorov}. It was
named due to the stochastic integral representation in terms of
standard Brownian motion, which was given by Mandelbrot and Van
Ness, \cite{bib:mandelbrot}.

Let $H$ be a constant in $(0,1).$ FBm $(B_{t}^{H})_{t\geq 0}$ with
Hurst parameter $H$ is a continuous and centered Gaussian process
with covariance function
$$E[B^{H}_{t}B^{H}_{s}]=\frac{1}{2}(t^{2H}+s^{2H}-|t-s|^{2H}).$$

For $H=1/2,$ fBm corresponds to a standard Brownian motion. A
standard fBm, $(B_{t}^{H})_{t\geq 0}$ has the following properties:
\begin{itemize}
  \item $B^{H}_{0}=0$ and $E[B^{H}_{t}]=0$ for all $t\geq 0.$
  \item $B^{H}$ has homogenous increments, that is
  $B^{H}_{t+s}-B^{H}_{s}$ has the same law as $B^{H}_{t}$ for $s,t\geq 0.$
  \item $B^{H}$ is a Gaussian process and $E[(B^{H}_{t})^{2}]=t^{2H}, t\geq
  0,$ for all $H\in (0,1).$
  \item $B^{H}$ has continuous trajectories.
\end{itemize}
\newpage
\textit{Stochastic integral representation}

Several stochastic integral representations have been developed for
the fBm. For example, it is proved that the following process is a
fBm with Hurst parameter $H\in (0,1),$
\begin{eqnarray}
B^{H}_{t}&=&\frac{1}{\Gamma(H+1/2)}\int_{\mathbb{R}}^{}((t-s)^{H-1/2}_{+}-(-s)^{H-1/2}_{+})dW_{s}\nonumber\\
&=&\frac{1}{\Gamma(H+1/2)}(\int_{-\infty}^{0}((t-s)^{H-1/2}-(-s)^{H-1/2})dW_{s}\nonumber\\
&+&\int_{0}^{t}(t-s)^{H-1/2}dW_{s})
\end{eqnarray}
where $W_{t}$ is a standard Brownian motion with $W_{0}=0$
considered on a probability space $(\Omega, \mathcal{F},
\mathbf{P}).$ The paths of $W$ are continuous, the increments of $W$
over disjoint intervals are independent Gaussian random variables
with zero-mean and with variance equal to the length of the
interval. And $\Gamma$ is the gamma function. Now let us omit the
constant $1/\Gamma(H+1/2)$ for simplicity and use the change of
variable $s=tu,$ then
\begin{eqnarray}E[(B^{H}_{t})^{2}]&=&\int_{\mathbb{R}}^{}((t-s)^{H-1/2}_{+}-(-s)^{H-1/2}_{+})^{2}d s\nonumber\\
&=&t^{2H}\int_{\mathbb{R}}^{}((1-u)^{H-1/2}_{+}-(-u)^{H-1/2}_{+})^{2}du\nonumber\\
&=&C(H)t^{2H}\end{eqnarray} Also
\begin{eqnarray}E[(B^{H}_{t}-B^{H}_{s})^{2}]&=&\int_{\mathbb{R}}^{}((t-u)^{H-1/2}_{+}-(s-u)^{H-1/2}_{+})^{2}d s\nonumber\\
&=&t^{2H}\int_{\mathbb{R}}^{}((t-s-u)^{H-1/2}_{+}-(-u)^{H-1/2}_{+})^{2}du\nonumber\\
&=&C(H)|t-s|^{2H}.\end{eqnarray} Hence
\begin{eqnarray}E[B^{H}_{t}B^{H}_{s}]&=&-\frac{1}{2}{E[(B^{H}_{t}-B^{H}_{s})^{2}]-E[(B^{H}_{t})^{2}]-E[(B^{H}_{s})^{2}]}\nonumber\\
&=&\frac{1}{2}(t^{2H}+s^{2H}-|t-s|^{2H}).\end{eqnarray}
\newpage
\textit{Correlation Between Two Increments}

For $H=\frac{1}{2},$ the process $(B^{H}_{t})_{t\geq0}$ corresponds
to a standard Brownian motion, in which the increments are
independent. For $H\neq 1/2$ the increments are not independent. By
the definition of fBm, we know the covariance between
$B^{H}(t+h)-B^{H}(t)$ and $B^{H}(s+h)-B^{H}(s)$ with $s+h\leq t$ and
$t-s=nh$ is
$$\rho_{H}(n)=\frac{1}{2}h^{2H}[(n+1)^{2H}+(n-1)^{2H}-2n^{2H}].$$

We observe that two increments of the form $B^{H}(t+h)-B^{H}(t)$ and
$B^{H}(t+2h)-B^{H}(t+h)$ are positively correlated for $H>1/2,$ and
they are negatively correlated for $H<1/2.$

\textit{Self Similarity Property}

Since the covariance function of fBm is homogenous of order $2H,$
fBm possess the self-similarity property, that is for any constant
$c>0,$
\begin{equation}\label{eq:scaling}(B^{H}_{ct})_{t\geq0}\stackrel{law}{=}(c^{H}B^{H}_{t})_{t\geq0}.\end{equation}
Note that, by taking $H=\frac{1}{2}$ we obtain the self similarity
property of standard Brownian motion $(W_{t})_{t\geq0}$ that is
\begin{equation}\label{eq:scalebm}(W_{at})_{t\geq0}\stackrel{law}{=}(a^{1/2}W^{H}_{t})_{t\geq0}.\end{equation}
Using equation $(\ref{eq:scalebm}),$ it was shown in
\cite{bib:yormakale} that the supremum of reflected standard
Brownian motion up to time $1$ is identical in law with the the
reciprocal of square root of hitting time of level $1.$ From there
on, the first moment of this supremum was calculated using the
properties of Normal distribution.

\newpage
\section{Main Results}

\textit{Summary of the main results}

In this section, we prove that the second moment of reflected fBm up
to time $a$ is bounded above by $a^{2H}.$ In our proof, as a new
approach we apply the properties of Gamma distribution to fBm
combined with the result given in \cite{bib:david} and use equation
(\ref{eq:scaling}).

Another new result given in this section is, we provide a lower
bound on the distribution of the supremum of fBm up to time $1$ and
up to fixed time $t.$ In this proof, we consider fBm up to a random
exponentially distributed time $T$ which is independent of the
process and using the self similarity property we obtain an upper
bound on the expected value of the supremum of fBm up to time $1$
and we combine this result with Markov's Inequality to find the
lower bound on the distribution of the supremum.

\textit{Notation}

Let $(B^{H}_{t})_{t\geq0}$ be a fBm defined on the probability space
$(\Omega, \mathcal{F}, \mathbf{P}).$ And, let us define the
reflected fBm around $0,$ that is denoted as
$(|B^{H}_{t}|)_{t\geq0},$
\begin{eqnarray}|B^{H}_{t}|_{t\geq0}:=\left\{\begin{array}{rl} B^{H}_{t}, ~~~~ B^{H}_{t}\geq 0\\ -B^{H}_{t}, ~~~~ B^{H}_{t}<0\end{array}\right\}\end{eqnarray}

Let $\tau_{1}:=\inf\{t\geq0 : |B^{H}_{t}|=1\}.$ In other words, $\tau_{1}$ is the first hitting time of level $1$ for the reflected fBm.

\begin{figure}[H]
\caption{Fractional Brownian motion and their reflections around
zero with Hurst parameters H=0.3,H=0.5,H=0.8}
\label{fig:fbmverfbmall}
\begin{center}
\includegraphics [type=eps,ext=.eps,read=.eps,width=15.5 cm,height=8.5 cm]{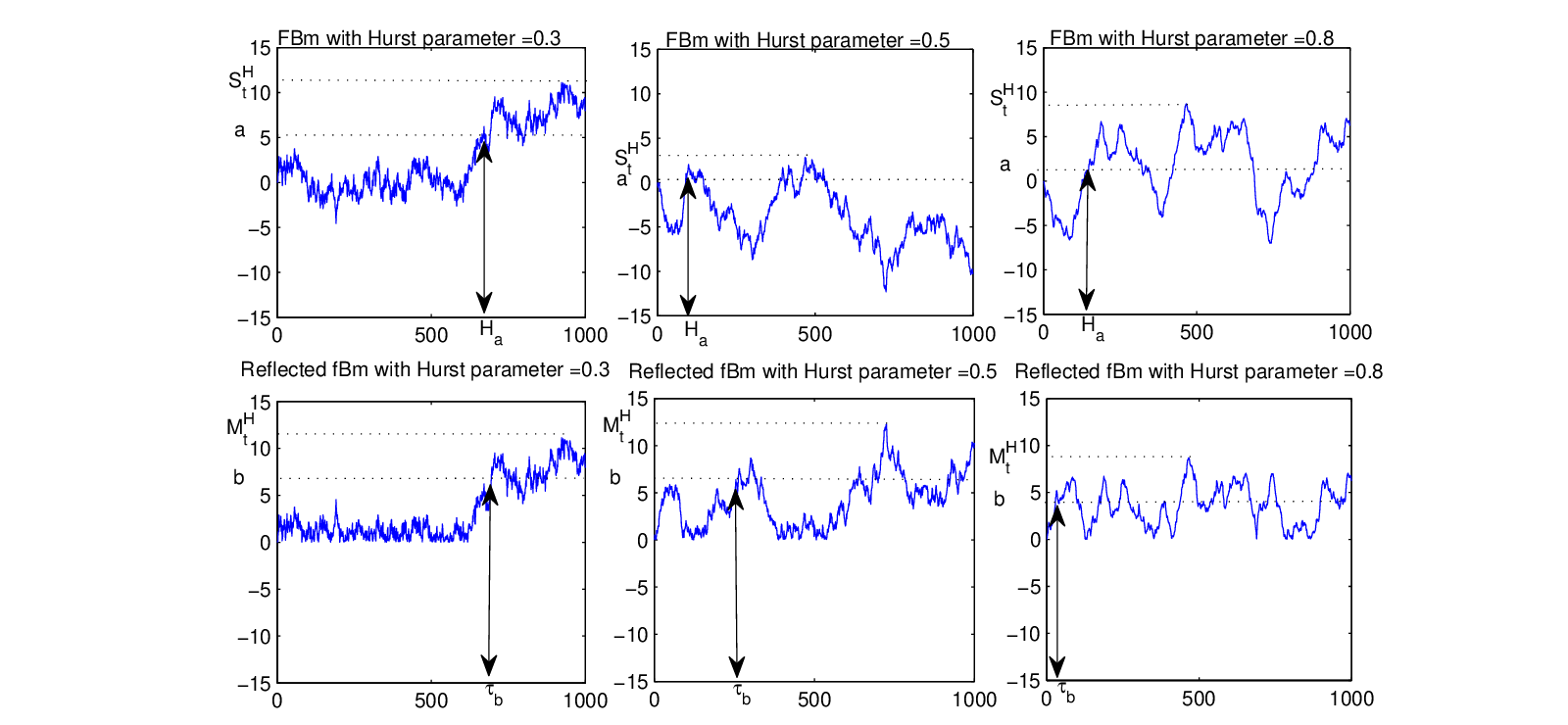}
\end{center}
\end{figure}

Let $H_{a}:=\inf\{t\geq0 : B^{H}_{t}=a\}$ be the first hitting time
of level $a$ where $a\in\mathbb{R_{+}}.$ Let $M^{H}_{t}:=\sup_{0\leq
v\leq t}|B^{H}_{v}|,$ that is the supremum of reflected fBm. And
similarly let $S^{H}_{t}:=\sup_{0\leq v\leq t}B^{H}_{v}$ be the
supremum of fBm.
\newpage
In Figure \ref{fig:fbmverfbmall}., we display sample paths of fBm
and reflected fBm with Hurst parameters $H=0.3,$ $H=0.5$ and $H=0.8$
from left to right. These sample paths are generated using the
Matlab code and the algorithm proposed by ABRY and SELLAN given in
\cite{bib:abrysellan}. This algorithm is a fast implementation of a
method proposed by SELLAN \cite{bib:sellan}, which is using sum
through a fractional wavelet basis. In this method, final paths
carry both the short-term and long-term correlation information. For
the details of simulations of fBm, one can also see a study by
CAGLAR, given in \cite{bib:caglar}.

\begin{theorem} For fBm with Hurst parameter $H>\frac{1}{2}$ and $a\in\mathbb{R_{+}}$
$$(M^{H}_{a})^{2}\stackrel{law}{=}(\frac{a}{\tau_{1}})^{2H}~~and~~(S^{H}_{a})^{2}\stackrel{law}{=}(\frac{a}{H_{1}})^{2H}$$
The second moment
$$E(M^{H}_{a})^{2}\leq a^{2H}$$
\end{theorem}
\textit{Proof.}

For $x>0,$ by equation (\ref{eq:scaling}) we obtain
\begin{eqnarray}\label{eq:law1}P((\frac{a}{\tau_{1}})^{2H}\leq x)&=&P(\tau_{1}\geq\frac{a}{x^{1/2H}})=P(\sup_{0\leq t\leq\frac{a}{x^{1/2H}}}|B^{H}_{t}|\leq 1)\\\nonumber
&=&P(\sup_{0\leq u\leq a}|B^{H}_{\frac{u}{x^{1/2H}}}|\leq
1)=P(\sup_{0\leq u\leq
a}|B^{H}_{u}|\leq\sqrt{x})=P((M^{H}_{a})^{2}\leq x)\end{eqnarray}
Hence, it is seen that
$$(M^{H}_{a})^{2}\stackrel{law}{=}(\frac{a}{\tau_{1}})^{2H}.$$

Now note that, \begin{eqnarray}\label{eq:law2}P(\sup_{0\leq u\leq
a}|B^{H}_{u}|\leq\sqrt{x})&=&P(\sup_{0\leq u\leq
a}B^{H}_{u}\leq\sqrt{x},\inf_{0\leq u\leq
a}B^{H}_{u}\geq\sqrt{x})\\\nonumber &\leq& P(\sup_{0\leq u\leq
a}B^{H}_{u}\leq\sqrt{x})=P((S^{H}_{a})^{2}\leq x)\end{eqnarray}

And by following the same argument given in (\ref{eq:law1}) it is
also observed that
$$(S^{H}_{a})^{2}\stackrel{law}{=}(\frac{a}{H_{1}})^{2H}.$$

Hence, combining the equations (\ref{eq:law1}) and (\ref{eq:law2})
the following inequality can be obtained
$$P((M^{H}_{a})^{2}\leq x)\leq P((S^{H}_{a})^{2}\leq x)=P((\frac{a}{H_{1}})^{2H}\leq x)$$

Now, by applying the properties of the Gamma distribution, for
$x>0,$ $k>0$ and $\theta>0$ we have
\begin{equation}\int^{\infty}_{0}\frac{x^{k-1}e^{-\frac{x}{\theta}}}{\Gamma(k)}dx=\theta^{k}.\end{equation}
Therefore
$$E((M^{H}_{a})^{2})\leq E((S^{H}_{a})^{2})=a^{2H}E((\frac{1}{H_{1}})^{2H})=a^{2H}\int^{\infty}_{0}E(e^{-xH_{1}^{2H}})dx$$
For a standard Brownian motion $W,$ it is well known that
$E(e^{-\lambda T_{a}})=e^{-a\sqrt{2\lambda}},$ where $T_{a}$ denotes
the first hitting time of level $a$ for standard Brownian motion. It
was shown in \cite{bib:david} that, for fBm with Hurst parameter
$H>\frac{1}{2}$ the following inequality holds
\begin{equation}\label{eq:bound}E(e^{-\lambda H^{2H}_{a}})\leq
e^{-a\sqrt{2\lambda}}\end{equation} for all $\lambda,a>0.$

Therefore, by equation (\ref{eq:bound})
\begin{equation}a^{2H}\int^{\infty}_{0}E(e^{-x\tau_{1}^{2H}})dx\leq
a^{2H}\int^{\infty}_{0}e^{-\sqrt{2x}}dx=a^{2H}\end{equation} As a
result, $E((M^{H}_{a})^{2})\leq a^{2H}.$ \hfill $\square$
\begin{corollary}For fBm with Hurst parameter $H>\frac{1}{2}$ and $a\in\mathbb{R_{+}}$
$$P(M^{H}_{a}>x)\leq P(S^{H}_{a}>x)\leq \frac{a^{H}}{x}$$\end{corollary}
\textit{Proof.}

Note that, by Jensen's inequality
$$E^{2}(S^{H}_{a})\leq E((S^{H}_{a})^{2})\leq a^{2H}$$
As a consequence, $E(M^{H}_{a})\leq E(S^{H}_{a})\leq a^{H}.$

Now, by Markov's inequality, it is seen that
$$P(M^{H}_{a}>x)\leq P(S^{H}_{a}>x)\leq \frac{a^{H}}{x}$$ \hfill $\square$

The bounds given above already are very useful information for the
investors because it tells the investor the probability of the
values of the supremum up to a fixed time.

In the next theorem a closer upper bound is obtained for the
distribution of the supremum of fBm which actually provides better
information for the investors.

\begin{theorem}\label{thm:bound} Consider fBm up to time $a,$ with Hurst parameter $H>\frac{1}{2},$
$$P(S^{H}_{a}\geq x)\leq \frac{\sqrt{2}a^{H}}{x\sqrt{\pi}}$$
\end{theorem}
\textit{Proof.}

Now, consider taking fBm up to time $T,$ where $T$ is exponentially
distributed random variable with mean $1/\lambda,$ independent of
underlying fBm, $B.$ Then using self similarity property of fBm,
\begin{equation}P(H^{2H}_{a}\leq T)=P(S_{T^{\frac{1}{2H}}}\geq a)=E(e^{-\lambda H^{2H}_{a}})\leq e^{-a\sqrt{2\lambda}}\end{equation}
and hence
\begin{equation}E(S_{T^{\frac{1}{2H}}})=\int_{0}^{\infty}P(S_{T^{\frac{1}{2H}}}\geq a)\leq \frac{1}{\sqrt{2\lambda}}\end{equation}
Again by the self similarity property we have,
\begin{equation}E(\sup_{0\leq u\leq T^{\frac{1}{2H}}}B_{u})=E(T^{1/2}\sup_{0\leq u\leq 1}B_{s})=E(\sqrt{T})E(S_{1})\end{equation}
where $E(T^{p/2})=\frac{\Gamma(\frac{2+p}{2})}{\lambda^{p/2}}.$

One can see \cite{bib:mainpaper}, for details related to taking a
random process up to a random time which is independent of the
underlying process.

Therefore we obtain
$$E(S_{1})\leq \frac{\sqrt{2}}{\sqrt{\pi}}$$
And by the Markov's Inequality for any $x>0,$ we can write
\begin{equation}P(S_{1}\geq x)\leq \frac{E(S_{1})}{x}=\frac{\sqrt{2}}{x\sqrt{\pi}}.\end{equation}
Now by using the scaling property, for any fixed time $a$ we observe
\begin{eqnarray}P(\sup_{0\leq u\leq
a}B^{H}_{u}\leq x)&=&P(a^{H}\sup_{0\leq u\leq 1}B^{H}_{u}\leq
x)=P(a^{H}S^{H}_{1}\leq x)\\ \nonumber &=&P(S^{H}_{1}\leq
\frac{x}{a^{H}})
\end{eqnarray}\hfill $\square$
\section{Application of the results to financial mathematics}

In finance, depending on the investor's demand, risk can be defined
in many ways. One of the commonly used definitions of risk is known
as $"downfalls"$ or $"maximum~~drawdown".$ Maximum drawdown is the
highest possible loss in the price of one share of the risky asset.
One of the fundamental results related to a standard Brownian motion
is L$\acute{e}$vy Theorem (see \cite{bib:yorkitap}). This theorem
leads us to the result that $"maximum~~drawdown",$ the highest
possible loss in the trajectories of the standard Brownian motion is
identical in law with the supremum of the reflected Brownian motion.
In other words, the $"maximum~~drawdown",$ before time $t$ is
defined as
$$D^{-}_{t}:=\sup_{0\leq u\leq v\leq t}(W_{u}-W_{v})=\sup_{0\leq v\leq t}(\sup_{0\leq u\leq v}W_{u}-W_{v})$$
and it satisfies $$(\sup_{0\leq u\leq
v}W_{u}-W_{v})\stackrel{law}{=}(\sup_{0\leq v \leq t}|W_{v}|)$$ by
the L$\acute{e}$vy isomorphism.

In the above result, we would like to replace standard Brownian
motion with fBm, because it is obvious that fBm is a more realistic
model than a standard Brownian motion for the risky assets. In the
literature, due to the complex nature of fBm there are not many
exact results on the distributions related to the trajectories. For
this reason, we also expect to obtain not exact results but some
bounds.

In the fBm notation let,
$$D^{-,H}_{t}:=\sup_{0\leq a'\leq a\leq t}(B_{a'}-B_{a})=\sup_{0\leq a\leq t}(\sup_{0\leq a'\leq a}B_{a'}-B_{a})$$
be the $"maximum~~drawdown".$

Our goal is to find the exact distribution or alternatively find
bounds for the distribution of $"maximum~~drawdown",$ which is
directly applicable to finance, however we have not been able to
obtain that, yet. In this study, we would like to present a lower
bound for the distribution of the difference between the supremum of
fBm and the value of fBm at time $a,$ which is the first step of
finding a lower bound for the distribution of $"maximum~~drawdown"$
of fBm up to time $t.$ And this following result will be extended to
finding the distribution of $"maximum~~drawdown"$ in the later
studies.

Now, let
\begin{eqnarray}\label{eq:difference}Y^{H}_{a}:=S^{H}_{a}-B^{H}_{a}\end{eqnarray}
be the difference between the supremum of fBm up to time $a$ and the
value of fBm at time $a.$

\begin{theorem} For $y\geq 0,$
$$P(Y^{H}_{a}\leq y)\geq 1-\frac{\sqrt{2}a^{H}}{y\sqrt{\pi}}$$
\end{theorem}
\textit{Proof.}

For $a\in\mathbb{R_{+}},$
\begin{eqnarray}
P(S^{H}_{a}\geq b)&=&P(S^{H}_{a}\geq b ,B^{H}_{a}\geq
b)+P(S^{H}_{a}\geq b ,B^{H}_{a}< b)\nonumber\\
&=& P(B^{H}_{a}\geq b)+P(S^{H}_{a}\geq b ,B^{H}_{a}< b)\leq
\frac{\sqrt{2}a^{H}}{x\sqrt{\pi}}
\end{eqnarray}
by Theorem (\ref{thm:bound}).

Thus we have, \begin{eqnarray}\label{eq:boundprob}0\leq
P(S^{H}_{a}\geq b ,B^{H}_{a}< b)\leq
\frac{\sqrt{2}a^{H}}{b\sqrt{\pi}}-P(B^{H}_{a}\geq b)\end{eqnarray}
For $b\geq x,$ $x\in\mathbb{R_{+}},$
\begin{eqnarray}&& P(B^{H}_{a}< x)-P(S^{H}_{a}< b ,B^{H}_{a}< x)\nonumber\\
&=&P(S^{H}_{a}\geq x ,B^{H}_{a}< x)\leq P(S^{H}_{a}\geq b
,B^{H}_{a}< b)\leq \frac{\sqrt{2}a^{H}}{b\sqrt{\pi}}-P(B^{H}_{a}\geq
b)\end{eqnarray} by (\ref{eq:boundprob}) and since $x\leq b.$
Therefore, we have
\begin{eqnarray}P(S^{H}_{a}< b ,B^{H}_{a}< x)\geq P(B^{H}_{a}<
x)+P(B^{H}_{a}\geq
b)-\frac{\sqrt{2}a^{H}}{b\sqrt{\pi}}\end{eqnarray} Now as a result,
we get
\begin{eqnarray}P(Y^{H}_{a}\leq y, B^{H}_{a}\leq x)\geq P(B^{H}_{a}<x)+P(B^{H}_{a}\geq y+x)-\frac{\sqrt{2}a^{H}}{(y+x)\sqrt{\pi}}\end{eqnarray}
by (\ref{eq:difference}), where $y=b-x$ and
\begin{eqnarray}P(Y^{H}_{a}\leq y,B^{H}_{a}\in dx)\geq
P(B^{H}_{a}\in dx)+\frac{dP(B^{H}_{a}\geq
y+x)}{dx}-\frac{\sqrt{2}a^{H}}{(y+x)^{2}\sqrt{\pi}}\nonumber\end{eqnarray}
finally
\begin{eqnarray}P(Y^{H}_{a}\geq y)\geq \int^{\infty}_{0}P(B^{H}_{a}\in
dx)+\int^{\infty}_{0}\frac{dP(B^{H}_{a}\geq
y+x)}{dx}dx-\int^{\infty}_{0}\frac{\sqrt{2}a^{H}}{(y+x)^{2}\sqrt{\pi}}dx\nonumber\end{eqnarray}
Hence,  for $y\geq0,$
\begin{eqnarray}
P(Y^{H}_{a}\leq y)&\geq&
\frac{1}{2}+\int^{\infty}_{0}dP(B^{H}_{a}\geq
y+x)-\int^{\infty}_{0}\frac{\sqrt{2}}{(y+x)^{2}\sqrt{\pi}}dx\nonumber\\
&=&\frac{1}{2}-\frac{\sqrt{2}a^{H}}{\sqrt{\pi}y}+\int^{\infty}_{0}dP(B^{H}_{a}\geq
y+x)
\end{eqnarray}
And by the Fundamental Theorem of Calculus,
\begin{eqnarray}P(Y^{H}_{a}\leq
y)\geq
\frac{1}{2}-\frac{\sqrt{2}a^{H}}{\sqrt{\pi}y}+\int^{\infty}_{0}\frac{1}{\sqrt{2\pi}}e^{\frac{-(x+y)^{2}}{2}}dx=
1-\frac{\sqrt{2}a^{H}}{\sqrt{\pi}y}\end{eqnarray}\hfill $\square$

\section{Conclusion}
The Black-Scholes model for the values of the log returns of a risky
stock using fBm is a more realistic model than using standard
Brownian motion, due the advantage of its capability of displaying
the dependence of increments. The investors of the risky stocks or
their derivatives, would naturally be interested in the highest
value of the asset, in a certain time period. Towards that end, we
have provided an upper bound on the second moment of the supremum of
reflected fBm as a first result, and based on this, we have given an
upper bound on the distribution of the supremum of reflected fBm.

Later, we have observed a lower bound on the distribution of the
supremum of fBm up to time $1$ as well as up to fixed time $a.$ As
an application of the second result, we have given the joint
distribution of the supremum of fBm up to time $a$ and the value of
fBm at time $a$ and obtained a lower bound on the distribution of
the difference between them. Obviously, this application is already
very useful for the investors. Based on the given results, we will
investigate the distribution of "maximum drawdown" because of its
important use as a measure of risk, in finance, as future work. Our
conjecture on this distribution is, that it is both related to the
distribution of the supremum of the reflected fBm and the
distribution of the difference between the supremum of fBm and the
value of fBm at the maturity time.

\section{Acknowledgement}
Ceren Vardar would like to thank Mine Caglar from Koc University,
Istanbul for her warm and generous support for this study. She would
also like to thank Craig L. Zirbel from Bowling Green State
University, Ohio and Gabor J. Szekely from National Science
Foundation, Washington D.C for introducing her to the subject.

\end{document}